\documentclass[a4paper,11pt]{amsart}
\usepackage{amssymb}
\usepackage{mathptmx}

\title[On a comparison of minimal log discrepancies]
{On a comparison of minimal log discrepancies \\ in terms of motivic integration}
\author{Masayuki Kawakita}

\address{Research Institute for Mathematical Sciences, Kyoto University, Kyoto 606-8502, Japan}
\email{masayuki@kurims.kyoto-u.ac.jp}

\newtheorem{theorem}{Theorem}[section]
\newtheorem{proposition}[theorem]{Proposition}
\newtheorem{lemma}[theorem]{Lemma}

\theoremstyle{definition}
\newtheorem{definition}[theorem]{Definition}
\newtheorem{example}[theorem]{Example}

\theoremstyle{remark}
\newtheorem{remark}[theorem]{Remark}

\numberwithin{equation}{section}

\newcommand{\bA}{{\mathbb A}}
\newcommand{\bG}{{\mathbb G}}
\newcommand{\bL}{{\mathbb L}}
\newcommand{\bQ}{{\mathbb Q}}
\newcommand{\bZ}{{\mathbb Z}}
\newcommand{\cB}{{\mathcal B}}
\newcommand{\cC}{{\mathcal C}}
\newcommand{\cD}{{\mathcal D}}
\newcommand{\cI}{{\mathcal I}}
\newcommand{\cJ}{{\mathcal J}}
\newcommand{\cM}{{\mathcal M}}
\newcommand{\cO}{{\mathcal O}}
\newcommand{\fm}{{\mathfrak m}}

\newcommand{\Hom}{{\mathcal H}om}
\newcommand{\Int}{\mathrm{Int}}
\newcommand{\Id}{\mathrm{Id}}
\newcommand{\mld}{\mathrm{mld}}
\newcommand{\Sch}{\mathrm{Sch}}
\newcommand{\chr}{\operatorname{char}}
\newcommand{\GL}{\operatorname{GL}}
\newcommand{\mult}{\operatorname{mult}}
\newcommand{\ord}{\operatorname{ord}}
\newcommand{\Spec}{\operatorname{Spec}}

\newcommand{\angl}[1]{\langle{#1}\rangle}

\begin{document}
\begin{abstract}
We formulate a comparison of minimal log discrepancies of a variety and its ambient space with appropriate boundaries in terms of motivic integration. It was obtained also by Ein and Musta\c{t}\v{a} independently.
\end{abstract}

\maketitle

\section{Introduction}
The purpose of this paper is to expand into an arbitrary variety the approach to the study of singularities in the minimal model program by means of motivic integration due to Ein, Musta\c{t}\v{a} and Yasuda. We introduce an ideal sheaf with rational exponent which measures how far a variety is from being local complete intersection. It naturally appears as a boundary in the adjunction of canonical divisors of an original variety and its ambient space. We prove the precise inversion of this adjunction, stated below with notation in Section \ref{sec:defect}.

\begin{theorem}\label{thm:main}
Let $X$ be a normal $\bQ$-Gorenstein closed subvariety of codimension $c$ of a smooth variety $A$, $\cI_X$ denote the corresponding ideal sheaf on $A$, and $\cD_X$ a weak l.c.i.\ defect $\bQ$-ideal sheaf of $X$. Let $Z$ be a closed proper subset of $X$. Then
\begin{align*}
\mld_Z(X,\cD_X)=\mld_Z(A,\cI_X^c).
\end{align*}
\end{theorem}

This theorem is a comparison of minimal log discrepancies.
Minimal log discrepancies play a pivotal role in the minimal model program. This program works with the provision of termination of flips thanks to a recent work \cite{HM05} of Hacon and McKernan, while the termination is reduced by Shokurov in \cite{S04} to two difficult conjectures on minimal log discrepancies. This is the main motivation to study a minimal log discrepancy, but also it is an interesting object itself because experience shows that it measures degree of singularities. For example, a normal surface is smooth if its minimal log discrepancy is greater than one, has Du Val singularities if it is at least one, and has quotient singularities if it is greater than zero. Shokurov conjectured its boundedness by dimension from above in \cite{S88}, but no bound has been obtained in general yet.

Ein, Musta\c{t}\v{a} and Yasuda applied motivic integration to the study of minimal log discrepancies in \cite{EM04}, \cite{EMY03}. They provided a description of minimal log discrepancies in terms of the space of arcs, and proved as applications precise inversion of adjunction and lower semi-continuity of minimal log discrepancies in the case of local complete intersection. For a normal $\bQ$-Gorenstein variety $X$ of dimension $d$ with $rK_X$ a Cartier divisor, the ideal sheaf $\cJ_{r,X}$ defined by the surjective map $(\Omega_X^d)^{\otimes r} \to \cJ_{r,X}\cO_X(rK_X)$ emerges naturally in the theory of motivic integration. Their method makes use of the fact that $\cJ_{1,X}$ coincides with the Jacobian ideal sheaf $\cJ'_X$, which is not difficult to calculate explicitly, provided that $X$ is local complete intersection. However in general the integral closure of $(\cJ'_X)^r$ is only contained in that of $\cJ_{r,X}$, of which the difference is to be studied in this paper.

We introduce the notion of weak l.c.i.\ defect $\bQ$-ideal sheaves $\cD_X$. A $\bQ$-ideal sheaf is just an extension of an ideal sheaf to that with a non-negative rational exponent. A $\cD_X$ has the property that its $r$-th power is represented by a usual ideal sheaf $\cD_{r,X}$ for which $\cJ_{r,X}\cD_{r,X}$ has the same integral closure as $(\cJ'_X)^r$ has. We show the existence of $\cD_X$ by a concrete construction, in the course of which we embed $X$ into various general l.c.i.\ schemes $Y = X \cup C^Y$ of the same dimension $d$. Then $\cD_{r,X}$ is achieved as the summation of ideal sheaves $\cO_X(-rC^Y|_X)$. In particular it is co-supported exactly on the non-l.c.i.\ locus of $X$, and one can say that $\cD_X$ measures how far $X$ is from being l.c.i. This $\cD_X$ naturally appears as a boundary to $X$ in the adjunction in Theorem \ref{thm:main}.

Theorem \ref{thm:main} is proved essentially along the idea of Ein, Musta\c{t}\v{a} and Yasuda. The new point is that we embed $X$ into a general l.c.i.\ scheme $Y = X \cup C^Y$ of the same dimension as in the construction of a weak l.c.i.\ defect $\bQ$-ideal sheaf. Then the embedding morphism from the space $J_\infty X$ of arcs on $X$ into that on $Y$ is a local isomorphism outside a subset of measure zero, by which the study of the space $J_\infty X$ is reduced to that of the space on $Y$. We avoid the detail analysis in \cite[Section 2]{EM04} of liftable jets thanks to Proposition \ref{prp:dimfibre}, so our proof simplifies that of the main theorem in \cite{EM04}.

Theorem \ref{thm:main} means that $X$ has better singularities by difference $\cD_X$ than the restriction of the pair $(A, \cI_X^c)$ has. It provides an affirmative answer to \cite[Conjecture 4.4]{T04}. Note that it is automatic to attach extra boundaries in the statement of Theorem \ref{thm:main}. On the other hand, S.\ Takagi's result \cite{T04} indicates that $X$ has worse singularities than the restriction of $(A, \cI_X)$ has. It is an interesting problem to find a boundary to $A$ between $\cI_X^c$ and $\cI_X$ to which the pair incidental corresponds to singularities of $X$, at least more approximately. This problem should guide us to deep investigation of minimal log discrepancies from the point of view of motivic integration. Such investigation is desirable in the sense that the result by means of motivic integration is at present the only plausible concrete evidence of precise inversion of adjunction and lower semi-continuity of minimal log discrepancies.

This paper is organised as follows. In Section \ref{sec:defect} we introduce and discuss weak l.c.i.\ defect $\bQ$-ideal sheaves. Section \ref{sec:motivic} is devoted to a quick review of motivic integration to state the characterisation Theorem \ref{thm:mld} of minimal log discrepancies by Ein, Musta\c{t}\v{a} and Yasuda. Our main Theorem \ref{thm:main} is proved in Section \ref{sec:comparison}. We work over an algebraically closed field $k$ of characteristic zero, although the discussions in Sections \ref{sec:defect} to \ref{sec:comparison} are valid even in positive characteristic with the provision of resolutions of singularities.

\section{Weak l.c.i.\ defect $\bQ$-ideal sheaves}\label{sec:defect}
The notion of $\bQ$-ideal sheaves on a scheme is an extension of ideal sheaves to those with non-negative rational exponents. A \textit{$\bQ$-ideal sheaf} is expressed as a finite product $\prod \cI_i^{a_i}$ of ideal sheaves with non-negative rational exponents $a_i$, and two expressions $\prod \cI_i^{a_i}$ and $\prod \cJ_j^{a_j}$ define the same $\bQ$-ideal sheaf if there exists a positive integer $r$ with any $ra_i, ra_j$ integral such that the $r$-th powers $\prod \cI_i^{ra_i}$ and $\prod \cJ_j^{ra_j}$ are the same as usual ideal sheaves. For a $\bQ$-ideal sheaf $\cI$, its \textit{denominator} is any positive integer $r$ such that the $r$-th power $\cI^r$ is expressed as a usual ideal sheaf, which we call a \textit{representative ideal sheaf} of $\cI$ to power of $r$. On a normal variety $X$ we identify an effective $\bQ$-Cartier divisor $D$ with a $\bQ$-ideal sheaf $(\cO_X(-rD))^{1/r}$, which does not equal $\cO_X(-D)$ in general, for a positive integer $r$ such that $rD$ is a Cartier divisor.

We introduce an equivalence relation in the category of $\bQ$-ideal sheaves on a fixed normal variety $X$, and write $\angl{\cI}$ for the equivalence class of a $\bQ$-ideal sheaf $\cI$. Two $\bQ$-ideal sheaves belong to the same equivalence class by definition if their representative ideal sheaves to power of a common denominator have the same integral closure. This equivalence relation is interpreted in terms of algebraic valuations. A valuation $v$ of the function field of $X$ is called an \textit{algebraic valuation} of $X$ if it is defined by a prime divisor $E$ on a smooth variety $\bar{X}$ equipped with a birational morphism $f \colon \bar{X} \to X$. Such a valuation is denoted by $v_E$, and the closure of $f(E)$ is called the \textit{centre} of $v_E$ on $X$ and denoted by $c_X(v_E)$. For a $\bQ$-ideal sheaf $\prod \cI_i^{a_i}$, its \textit{multiplicity} along $v_E$ is defined as $\sum a_i \mult_E \cI_i$, which is independent of expressions. Then two $\bQ$-ideal sheaves are equivalent if and only if their multiplicities along any algebraic valuation coincide. In particular multiplicities along algebraic valuations are well-defined in the set of equivalence classes of $\bQ$-ideal sheaves. We define a partial order $\subset$ on this set so that $\angl{\cI} \subset \angl{\cJ}$ if $\mult_E \cI \ge \mult_E \cJ$ for every algebraic valuation $v_E$.

In the minimal model program, it is natural to consider a pair $(X, \cI)$ of a normal $\bQ$-Gorenstein variety $X$ and a $\bQ$-ideal sheaf $\cI$ on $X$. Let $v_E$ be an algebraic valuation defined by a prime divisor $E$ on $\bar{X}$ equipped with $f \colon \bar{X} \to X$, and write $K_{\bar{X}} = f^*K_X+A$ with an exceptional $\bQ$-divisor $A$. Then the \textit{log discrepancy} $a_E(X,\cI)$ of $v_E$ with respect to $(X, \cI)$ is defined as $1 + \mult_E A - \mult_E \cI$. For a closed subset $Z$ of $X$, we define the \textit{minimal log discrepancy} $\mld_Z(X,\cI)$ of $(X, \cI)$ over $Z$ as the infimum of $a_E(X,\cI)$ for $c_X(v_E) \subset Z$. Minimal log discrepancies take values in $\bQ_{\ge0} \cup \{-\infty\}$when $X$ has dimension at least two, but for convenience we set $\mld_Z(X,\cI):=-\infty$ even in dimension one if the infimum is negative. We say that the pair $(X, \cI)$ is \textit{log canonical} near $Z$ if $\mld_U(U,\cI) \ge 0$ for a neighbourhood $U$ of $Z$, which is equivalent to $\mld_Z(X,\cI) \ge 0$. Similarly we say that $(X, \cI)$ is \textit{kawamata log terminal} near $Z$ if $\mld_U(U,\cI) >0$, which is stronger than $\mld_Z(X,\cI) >0$.

Let $X$ be a normal $\bQ$-Gorenstein variety of dimension $d$, and $r$ any positive integer such that $rK_X$ is a Cartier divisor. In application of motivic integration to the study of singularities in the minimal model program, the ideal sheaf $\cJ_{r,X}$ defined below plays a pivotal role. Consider the natural map $(\Omega_X^d)^{\otimes r} \to \cO_X(rK_X)$ and define$\cJ_{r,X}$ so that the image of this map is $\cJ_{r,X}\cO_X(rK_X)$. We denote by $\cJ_X$ the $\bQ$-ideal sheaf $(\cJ_{r,X})^{1/r}$, which is independent of $r$. If $X$ is local complete intersection (l.c.i.\ for short), then $X$ is Gorenstein and the sheaf $\cJ_{1,X}$ coincides with the \textit{Jacobian ideal sheaf} $\cJ'_X$, which is the $d$-th Fitting ideal sheaf of $\Omega_X^1$. However, in general only the inclusion $\angl{\cJ'_X} \subset \angl{\cJ_X}$ holds and this inclusion is strict exactly on the non-l.c.i.\ locus of $X$, as we will see below.

\begin{definition}
A \textit{weak l.c.i.\ defect $\bQ$-ideal sheaf} $\cD_X$ of a normal $\bQ$-Gorenstein variety $X$ is any $\bQ$-ideal sheaf which satisfies the equality $\angl{\cJ'_X}=\angl{\cJ_X}\angl{\cD_X}$.
\end{definition}

We do not know if such a $\bQ$-ideal sheaf $\cD_X$ exists a priori. We will construct $\cD_X$ concretely on the germ $x \in X$ at a closed point $x$ in a form applicable to the proof of our main theorem. Fix an ambient space $A$ of $X$ which is smooth, and for each closed subscheme $V$ of $A$ we denote by $\cI_V$ the corresponding ideal sheaf on $A$. We take generally a subscheme $Y$ of $A$ which contains $X$ and is l.c.i.\ of dimension $d$. Then by Bertini's theorem $Y$ is the scheme-theoretic union of $X$ and another variety, say $C^Y$, of dimension $d$, whence $\cI_Y=\cI_X \cap \cI_{C^Y}$ and $\cI_{C^Y}\cO_X$ equals the conductor ideal sheaf $\cC_{X/Y}:=\Hom_{\cO_Y}(\cO_X, \cO_Y)$. By Grothendieck duality we have
\begin{align}\label{eqn:Gduality}
\omega_X=\cI_{C^Y}\omega_Y\cO_X.
\end{align}
In particular the subscheme $D^Y:=C^Y|_X$ on $X$ which corresponds to the ideal sheaf $\cI_{C^Y}\cO_X$ is a $\bQ$-Cartier divisor, and
\begin{align}\label{eqn:rGduality}
\cO_X(rK_X)=\cO_X(-rD^Y)\omega_Y^{\otimes r}.
\end{align}
Take the surjective maps
\begin{align*}
(\Omega_A^d)^{\otimes r} \twoheadrightarrow \cJ_{r,Y}\omega_Y^{\otimes r} \twoheadrightarrow \cJ_{r,X}\cO_X(rK_X).
\end{align*}
Here the ideal sheaf $\cJ_{r,Y}$ is defined also for $Y$ by the surjective map $(\Omega_Y^d)^{\otimes r} \twoheadrightarrow \cJ_{r,Y}\omega_Y^{\otimes r}$, and in this case $\cJ_{r,Y}$ equals the $r$-th power $(\cJ'_Y)^r$ of the Jacobian ideal sheaf $\cJ'_Y$. Hence by (\ref{eqn:rGduality}) we have
\begin{align}\label{eqn:descriptJ}
(\cJ'_Y\cO_X)^r=\cJ_{r,X} \cdot \cO_X(-rD^Y).
\end{align}
Now take summation of equalities (\ref{eqn:descriptJ}) for all the general l.c.i.\ $Y$, so 
\begin{align}\label{eqn:sum_descriptJ}
\sum_Y (\cJ'_Y\cO_X)^r=\cJ_{r,X} \cdot \sum_Y \cO_X(-rD^Y).
\end{align}
The equivalence class $\angl{\sum_Y (\cJ'_Y\cO_X)^r}$ of its left-hand side equals $\angl{(\cJ'_X)^r}$ by the relation $\cJ'_X=\sum_Y \cJ'_Y\cO_X$, whence we obtain that
\begin{align*}
\angl{(\cJ'_X)^r}=\angl{\cJ_X^r}\angl{\sum_Y \cO_X(-rD^Y)}.
\end{align*}
Therefore the $\bQ$-ideal sheaf $(\sum_Y \cO_X(-rD^Y))^{1/r}$ is a weak l.c.i.\ defect $\bQ$-ideal sheaf of $X$.

We have $\sum_Y \cO_X(-rD^Y) \subset [(\cJ'_X)^r \colon \cJ_{r,X}]$ by (\ref{eqn:sum_descriptJ}), which implies that $\angl{\cD_X} \subset \angl{[(\cJ'_X)^r \colon \cJ_{r,X}]^{1/r}}$. On the other hand the inclusion $(\cJ'_X)^r \supset [(\cJ'_X)^r \colon \cJ_{r,X}] \cJ_{r,X}$ implies that $\angl{\cD_X} \supset \angl{[(\cJ'_X)^r \colon \cJ_{r,X}]^{1/r}}$. Hence $\angl{\cD_X}=\angl{[(\cJ'_X)^r \colon \cJ_{r,X}]^{1/r}}$ so we may choose $\cD_X=[(\cJ'_X)^r \colon \cJ_{r,X}]^{1/r}$, which can be globalised.

The discussions above are summarised as follows.

\begin{proposition}
A weak l.c.i.\ defect $\bQ$-ideal sheaf $\cD_X$ exists, and it is co-supported exactly on the non-l.c.i.\ locus of $X$. 
\end{proposition}

Presumably the ideal sheaf $\cD'_X:=\sum_Y \cC_{X/Y}$ is a more natural sheaf which measures how far $X$ is from being l.c.i., because it is defined for an arbitrary variety. We call $\cD'_X$ the \textit{l.c.i.\ defect ideal sheaf}. $\cD'_X$ is also a weak l.c.i.\ defect $\bQ$-ideal sheaf when $X$ is Gorenstein, but in general only the inclusion $\angl{\cD'_X} \subset \angl{\cD_X}$ holds. The difference of them stems from that of $\cO_X(K_X)^{\otimes r}$ and $\cO_X(rK_X)$.

Below we compute $\cD_X$ in two examples of non-l.c.i.\ Gorenstein singularities. As Example \ref{exl:1/3} shows, $\cD_X$ is in general too high to derive benefit from Theorem \ref{thm:main} if the codimension $c$ is large, whereas Example \ref{exl:codim3} supports that $\cD_X$ is relatively moderate when $c$ is small.

\begin{example}\label{exl:1/3}
Let $o \in X = \frac{1}{3}(1,1,1)$ be the quotient of $\bA^3 = \Spec k[x_1,x_2,x_3]$ by the action of $\bZ/(3)$ given by $x_i \mapsto \zeta x_i$, where $\zeta$ is a primitive cubic root of unity ($\chr k \neq 3$ assumed). Then $X$ is Gorenstein and admits an embedding $o \in X \subset A=\bA^{10}$ via $\bZ/(3)$-invariant monomials of degree $3$ in $x_1,x_2,x_3$. By $\omega_X = \cO_Xdx_1 \wedge dx_2 \wedge dx_3$ we have $\cJ_{1,X}=\fm^2$, where $\fm$ is the maximal ideal sheaf of the origin. On the other hand one can easily check that $((x_1^3)^7,(x_2^3)^7,(x_3^3)^7) \subset \cJ'_X \subset \fm^7$, whence $\cJ'_X$ has integral closure $\fm^7$. Therefore $\cD_X=\fm^5$ is a weak l.c.i.\ defect $\bQ$-ideal sheaf of $X$.
\end{example}

\begin{example}\label{exl:codim3}
A non-l.c.i.\ Gorenstein singularity has codimension at least $3$ to its ambient space, and such a singularity $x \in X$ of codimension $3$ is explicitly described by means of pfaffians thanks to Buchsbaum and Eisenbud in \cite{BE77}. For an ambient space $A$ in which $X$ has codimension $3$, there exists an alternating map $f \colon \cO_A^n \to \cO_A^n$ with $n$ odd such that $X$ is given by the ideal sheaf generated by all the $(n-1)$-st order pfaffians of $f$. We claim that the l.c.i.\ defect ideal sheaf $\cD'_X$ is generated by all the $(n-3)$-rd order pfaffians of $f$. For an $n$-th alternating matrix $M$, we let $p(M)$ denote the pfaffian of $M$ and $M_{i_1 \ldots i_j}$ the minor of $M$ obtained by deletion of $i_1,\ldots,i_j$-th rows and columns. The incidental alternating matrix $Q^M=(q_{ij})$ with $q_{ij} := (-1)^{i+j}p_{ij}(M)$ for $i<j$ satisfies that $Q^MM = MQ^M = p(M)\Id_n$. Take indeterminates $m_{ij}$ for $1 \le i < j \le n$ and make an $n$-th alternating matrix $M=(m_{ij})$ by setting $m_{ji}:=-m_{ij}$ and $m_{ii}:=0$. The claim is reduced to the generic case in which $A$ is the affine space $\bA_R^{(n-1)n/2} = \Spec R[\underline{m}]$ over a regular ring $R$ finitely generated over $k$ and $X$ is given by the ideal sheaf $\cI_X=(p_1,\ldots,p_n)$, where we set $p_{i_1 \ldots i_j}:=p(M_{i_1 \ldots i_j})$. Consider the complete intersection subscheme $Y$ of $A$ defined by the ideal sheaf $\cI_Y=(p_1,p_2,p_3)$. It suffices to prove that the conductor ideal sheaf $\cC_{X/Y}$, principal by (\ref{eqn:Gduality}), equals $p_{123}\cO_X$. By $M\vec{p}=0$ with $\vec{p} := {}^t(-p_1, \ldots, (-1)^np_n)$ we have
$0=
\begin{pmatrix}
0 & Q^{M_{123}}
\end{pmatrix}
M\vec{p}=
\begin{pmatrix}
\ast & p_{123}\Id_{n-3}
\end{pmatrix}
\vec{p}$,
whence $p_{123} \in \cC_{X/Y}$. Now one can prove the equality $\cC_{X/Y}=p_{123}\cO_X$ by induction on $n$, for various $i_0,j_0$ by cutting objects out with hyperplanes $m_{ij}=0$ with $\emptyset \neq \{i,j\} \cap \{i_0,j_0\} \neq \{i_0,j_0\}$ and restricting to the locus $m_{i_0j_0}\neq0$.
\end{example}

\section{Review of motivic integration}\label{sec:motivic}
Motivic integration was invented by Kontsevich on a smooth variety in \cite{Kn95} as an analogue of $p$-adic integration, and was developed by Denef and Loeser on an arbitrary variety in \cite{DL99}. It is an integration on the space of arcs which takes value in an extension of the Grothendieck ring.

Let $\Sch/k$ denote the category of schemes of finite type over $\Spec k$. For a non-negative integer $n$, the covariant functor $\Sch/k \to \Sch/k$ which sends $X$ to $X \times_{\Spec k} \Spec k[t]/(t^{n+1})$ has the right adjoint functor $J_n$. A $k$-valued point of the scheme $J_nX$ corresponds to a morphism $\Spec k[t]/(t^{n+1}) \to X$, which we call a \textit{jet} of order $n$ on $X$. $J_nX$ admits a natural $\bA^1$-action. By $a \in k$ a jet $\gamma_n \colon \Spec k[t]/(t^{n+1}) \to X$ is sent to the jet $\gamma_n \circ g_a$ where $g_a \colon \Spec k[t]/(t^{n+1}) \to \Spec k[t]/(t^{n+1})$ is given by $t \mapsto at$. There exists a natural morphism $J_{n+1}X \to J_nX$, and we obtain a scheme $J_\infty X$, not of finite type in general, by taking the inverse limit of these $J_nX$. A $k$-valued point of $J_\infty X$ corresponds to a morphism $\Spec k[[t]] \to X$, which we call an \textit{arc} on $X$. $J_\infty X$ also admits a natural $\bA^1$-action. We set $\pi_n^X \colon J_\infty X \to J_nX$ and $\pi_{nm}^X \colon J_mX \to J_nX$ for $m \ge n$.

We have introduced the space on which motivic integration is to be defined, so we then define the space in which it takes value. The \textit{Grothendieck ring} $K_0(\Sch/k)$ is the quotient group of the free abelian group generated by isomorphism classes $[X]$ of $\Sch/k$ with relations $[X]=[Y]+[X \setminus Y]$ whenever $Y$ is a closed subscheme of $X$, endowed with the ring structure defined by $[X][Y]:=[X \times_{\Spec k} Y]$. Set $\bL:=[\bA^1]$ and take the localisation $\cM$ of $K_0(\Sch/k)$ by the multiplicatively closed set $\{\bL^n\}_{n \ge 0}$. We extend $\cM$ to the ring $\cM_\bQ:= \bigoplus_{q\in \bQ \cap [0,1)}\cM\bL^q$ in a natural way. The ring $\cM_\bQ$ has a descending filtration $F_n\cM_\bQ$ for $n \in \bQ$ where $F_n\cM_\bQ$ is the subgroup generated by all the $[X]\bL^q$ with $\dim X +q \le -n$. This filtration satisfies that $F_m\cM_\bQ \cdot F_n\cM_\bQ \subset F_{m+n}\cM_\bQ$, whence the completion $\widehat{\cM}_\bQ$ of $\cM_\bQ$ by this filtration becomes a topological ring equipped with an induced filtration $F_n\widehat{\cM}_\bQ$. The \textit{dimension} of an element of $\widehat{\cM}_\bQ$ is defined as the infimum of $n$ such that it is contained in $F_{-n}\widehat{\cM}_\bQ$.

We will define measurable sets and their motivic measures on the space of arcs $J_\infty X$ on a scheme $X$ of finite type of dimension $d$. We say that a subset $S$ of $J_\infty X$ is \textit{stable} at level $n$ if
\begin{enumerate}
\item
$\pi_n^X(S)$ is a constructible subset of $J_nX$,
\item
$S=(\pi_n^X)^{-1}(\pi_n^X(S))$, and
\item
for any $m \ge n$ the projection $\pi_{m+1}^X(S) \to \pi_m^X(S)$ is a piecewise trivial fibration with fibres $\bA^d$, where a fibration is said to be \textit{piecewise trivial} if it becomes trivial after a stratification of its base.
\end{enumerate}
The family $\cB_{0,X}$ of stable subsets of $J_\infty X$ is closed under finite intersections and unions. We define a function $\mu_{0,X} \colon \cB_{0,X} \to \widehat{\cM}_\bQ$ by sending $S$ to $[\pi_n^X(S)]\bL^{-(n+1)d}$. We say that a subset $S$ of $J_\infty X$ is \textit{measurable} if there exist subsets $S_i$ of $S$, stable subsets $T_i$ of $J_\infty X$ and subvarieties $Y_{ij}$ of $X$ with $\dim Y_{ij} <d$ for $i,j \ge 1$ such that
\begin{enumerate}
\item
$S$ is the disjoint union $\bigsqcup_{i \ge 1} S_i$,
\item
$S_i \circleddash T_i \subset \bigcup_{j \ge 1} J_\infty Y_{ij}$, where $S \circleddash T$ denotes the symmetric difference $(S \setminus T) \sqcup (T \setminus S)$ of $S$ and $T$, and
\item
the summation $\sum_{i \ge 1} \mu_{0,X}(T_i)$ has a limit.
\end{enumerate}
Then we define $\mu_X(S)$ as this limit, which is independent of $S_i$, $T_i$ and $Y_{ij}$. The family $\cB_X$ of measurable subsets of $J_\infty X$ forms a finitely additive class. We call the function $\mu_X \colon \cB_X \to \widehat{\cM}_\bQ$ the \textit{motivic measure}. Here we have extended $\cB_{0,X}$ to $\cB_X$ so that a subset $S$ of $J_\infty X$ has measure zero if and only if $S \subset \bigcup_{i \ge 1} J_\infty Y_i$ for subvarieties $Y_i$ with $\dim Y_i <d$.

Let $S$ be a subset of $J_\infty X$. A function $\alpha \colon S \to \bQ \cup \{\pm\infty\}$ is called a \textit{measurable function} if its fibres are measurable and $\mu_X(\alpha^{-1}(\infty))=\mu_X(\alpha^{-1}(-\infty))=0$. For such $\alpha$ we formally define the \textit{motivic integration} of $\bL^\alpha$ by
\begin{align*}
\int_S\bL^\alpha d\mu_X:=\sum_{q \in \bQ}\mu_X(\alpha^{-1}(q))\bL^q.
\end{align*}
We say that $\bL^\alpha$ is \textit{integrable} if this summation has a limit.
 
Let $\cI$ be a $\bQ$-ideal sheaf on $X$ co-supported on a locus of dimension less than $d$, with a representative ideal sheaf $\cI_{r}$ to power of a denominator $r$. The \textit{order} function $\ord_\cI \colon J_\infty X \to \frac{1}{r}\bZ_{\ge 0} \cup \{\infty\}$ is defined by the relation $\gamma^{-1}\cI_{r} \cdot \cO_{\Spec k[[t]]}=(t^{r\ord_\cI(\gamma)})$ for each arc $\gamma \colon \Spec k[[t]] \to X$, where we set $t^\infty:=0$. Then $\ord_\cI$ is a measurable function.

The most useful property of motivic integration is the formula of transformation. Let $f \colon \bar{X} \to X$ be a proper birational morphism from a smooth variety $\bar{X}$ to a variety $X$ of dimension $d$, and $f_\infty \colon J_\infty \bar{X} \to J_\infty X$ the induced morphism. Consider a natural map $f^* \Omega_X^d \to \omega_{\bar{X}}$ and define an ideal sheaf $\cJ_f$ on $\bar{X}$ so that the image of this map is $\cJ_f\omega_{\bar{X}}$. Let $S$ be a subset of $J_\infty X $ and $\alpha$ a measurable function on $S$. Then the function $\alpha \circ f_\infty$ on $f_\infty^{-1}(S)$ becomes a measurable function. The following \textit{formula of transformation} associates the integration of $\bL^\alpha$ to an integration on $f_\infty^{-1}(S)$, including integrability.

\begin{theorem}\label{thm:transformation}
\begin{align*}
\int_S\bL^\alpha d\mu_X=
\int_{f_\infty^{-1}(S)}\bL^{\alpha \circ f_\infty-\ord_{\cJ_f}}d\mu_{\bar{X}}.
\end{align*}
\end{theorem}

Ein, Musta\c{t}\v{a} and Yasuda applied motivic integration to the study of minimal log discrepancies in \cite{EMY03} by providing a description of minimal log discrepancies in terms of the space of arcs.

\begin{theorem}\label{thm:mld}
Let $(X, \cI)$ be a pair of a normal $\bQ$-Gorenstein variety $X$ and a $\bQ$-ideal sheaf $\cI$ on $X$, and $Z$ a closed subset of $X$. For a non-negative rational number $a$, the following are equivalent.
\begin{enumerate}
\item
$\mld_Z(X,\cI)<a$.
\item
There exists a stable subset $S$ of $J_\infty X$ mapped into $Z$ such that $\ord_{\cJ_X\cI}$ is constant on $S$ and such that $\dim \mu_X(S)+\ord_{\cJ_X\cI}S > -a$.
\end{enumerate}
In particular $\mld_Z(X,\cI)$ is an invariant on the formal scheme of $X$ along $Z$.
\end{theorem}

Their original description uses codimension of $S$, but it would be more natural to adopt dimension in the sense of motivic integration.

\begin{remark}
Theorem \ref{thm:mld} is concluded from a direct description of $\mld_Z(X,\cI)$ using the motivic integrations
\begin{align*}
\Int(\epsilon):=\int_{(\pi_0^X)^{-1}(Z)} \bL^{(1-\epsilon) \ord_{\cJ_X\cI}} d\mu_X.
\end{align*}
The pair $(X, \cI)$ is log canonical near $Z$ if and only if $\Int(\epsilon)$ has a limit for any positive rational number $\epsilon$, in which case $\mld_Z(X,\cI) = - \lim_{\epsilon \to +0} \dim \Int(\epsilon)$. If $(X, \cI)$ is kawamata log terminal near $Z$, then $\Int(0)$ has a limit and $\mld_Z(X,\cI) = - \dim \Int(0)$.
\end{remark}

\section{Comparison of minimal log discrepancies}\label{sec:comparison}
The purpose of this section is to prove Theorem \ref{thm:main}. The statement is local so we discuss on the germ $x \in Z \subset X \subset A$ at a closed point $x$. Denote by $d$ the dimension of $X$ and fix a positive integer $r$ such that $rK_X$ is a Cartier divisor. As in the construction of a weak l.c.i.\ defect $\bQ$-ideal sheaf $\cD_X$ in Section \ref{sec:defect}, we take generally a subscheme $Y$ of $A$ which contains $X$ and is l.c.i.\ of dimension $d$, so $Y$ is the scheme-theoretic union of $X$ and another variety $C^Y$, and $D^Y:=C^Y|_X$ is a $\bQ$-Cartier divisor on $X$. We denote by $\cI_V$ the ideal sheaf on $A$ which corresponds to a closed subscheme $V$ of $A$.

Take $c$ general prime divisors $H_1^Y,\ldots,H_c^Y$ on $A$ which contain $Y$, so $Y$ is the complete intersection of them. The Grothendieck duality (\ref{eqn:Gduality}) provides a natural adjunction
\begin{align}\label{eqn:duality}
K_X+D^Y=(K_A+\sum_{1 \le i \le c}H^Y_i)|_X.
\end{align}
By general choice of $Y$ and $H_i^Y$, we have
\begin{align*}
\mld_Z(X,\cD_X)   &= \mld_Z(X,D^Y), \\
\mld_Z(A,\cI_X^c) &= \mld_Z(A,\cI_Y^c) = \mld_Z(A,\sum_{1 \le i \le c}H_i^Y).
\end{align*}
Take a resolution $f \colon \bar{A} \to A$ such that
\begin{enumerate}
\item
its exceptional locus on $\bar{A}$ is a divisor $\sum_jE_j$,
\item
$\sum_i\bar{H}_i^Y + \sum_jE_j$ is a simple normal crossing divisor, where $\bar{H}_i^Y$ is the birational transform of $H_i^Y$,
\item
$\bigcap_{1 \le i \le c}\bar{H}_i^Y$ is the disjoint union of smooth varieties $\bar{X}$ and $\bar{C^Y}$ which are resolutions of $X$ and $C^Y$, and
\item
$E_j \cap \bar{X} \neq \emptyset$ implies that $f(E_j) = f(E_j|_{\bar{X}})$.
\end{enumerate}
Write $K_{\bar{A}} + \sum_i\bar{H}^Y_i + \sum_jE_j = f^*(K_A+\sum_iH_i^Y) + \sum_ja_jE_j$ with $a_j:=a_{E_j}(A,\sum_iH_i^Y)$. Then $K_{\bar{X}} + \sum_jE_j|_{\bar{X}} = f^*(K_X+D^Y) + \sum_ja_jE_j|_{\bar{X}}$ by (\ref{eqn:duality}). In particular the inequality $\mld_Z(X,D^Y) \ge \mld_Z(A,\sum_iH_i^Y)$ is trivial. It suffices to prove the converse inequality.

Suppose that $\mld_Z(A,\sum_iH_i^Y) < a$ for fixed $a \in \bQ_{\ge0}$. Our goal is to prove the inequality $\mld_Z(X,D^Y)<a$. We adopt the approach by means of motivic integration due to Ein, Musta\c{t}\v{a} and Yasuda. By the characterisation Theorem \ref{thm:mld} of minimal log discrepancies, there exist a stable subset $S_\infty^o$ of $J_\infty A$ at level $n$ mapped into $Z$ and an integer $p \le n$ such that $\ord_{H_i}$ takes constant value $p$ on $S_\infty^o$ for every $i$, and for any $m \ge n$
\begin{align}\label{eqn:dimS}
\dim S_m^o-(m+1)(d+c)+cp > -a,
\end{align}
where we set $S_m^o:=\pi_m^A(S_\infty^o)$. We may assume that $S_n^o$ is a $\bG_m$-invariant irreducible subset of $J_n A$. Let $S_\infty$ denote the closure of $S_\infty^o$ in $J_\infty A$ and $S_m$ that of $S_m^o$ in $J_mA$. Then $S_\infty$ is also stable at level $n$ by the smoothness of $A$. Let $W$ denote the singular locus of $Y$ and set $T'_\infty:=S_\infty \cap J_\infty X \setminus J_\infty W$.

\begin{lemma}
$T'_\infty \neq \emptyset$.
\end{lemma}

\begin{proof}
We follow the argument in the proof of \cite[Lemma 3.2]{EMY03}. Take a resolution $f \colon \bar{X} \to X$ and set $f_\infty \colon J_\infty \bar{X} \to J_\infty X$. $\pi_0^{\bar{X}}$ admits a natural zero section $s \colon \bar{X}=J_0 \bar{X} \to J_\infty \bar{X}$. $f_\infty^{-1}(S_\infty \cap J_\infty X)$ is $\bG_m$-invariant and closed, so $f_\infty^{-1}(S_\infty \cap J_\infty X)$ admits a natural $\bA^1$-action and $f_\infty^{-1}(S_\infty \cap J_\infty X)$ is sent by $0$ to $s(f^{-1}(\pi_0^A(S_\infty)))$. Hence $f_\infty^{-1}(S_\infty \cap J_\infty X) \neq \emptyset$. Since $\bar{X}$ is smooth, $f_\infty^{-1}(S_\infty \cap J_\infty X)$ is stable at level $n$ and $\mu_{\bar{X}}(f_\infty^{-1}(S_\infty \cap J_\infty X)) \neq 0$. Thus $\mu_X(S_\infty \cap J_\infty X) \neq 0$ by Theorem \ref{thm:transformation} and the lemma follows.
\end{proof}

Let $e$ denote the minimum value of $\ord_{\cJ'_Y}$ on $T'_\infty$, which is not infinity because $\cJ'_Y$ is co-supported on $W$. Set $T_\infty^o := T'_\infty \cap (\ord_{\cJ'_Y})^{-1}(e)$ and $T_m^o := \pi_m^X(T_\infty^o)$. We may assume that $n \ge re$ by replacing $n$ if necessary. For $m \ge n$ we set $U_m^o:=S_m \cap J_mY \cap (\ord_{\cJ'_Y})^{-1}(\bZ_{\le e}) \cap (\ord_{\cI_{C^Y}})^{-1}(\bZ_{\le re})$, which is an open subset of $S_m \cap J_mY$. Note that $T_m^o \subset U_m^o$ by (\ref{eqn:descriptJ}) and $\cI_{rD^Y} \subset \cI_{D^Y} = \cI_X+\cI_{C^Y}$.

\begin{lemma}\label{lem:stable}
$T_\infty^o$ is a stable subset of both $J_\infty X$ and $J_\infty Y$ at level $n$.
\end{lemma}

\begin{proof}
Take any $m \ge n$. We have the inclusion $(\pi_m^Y)^{-1}(T_m^o) \subset J_\infty Y \setminus J_\infty C^Y$. Since $J_\infty Y=J_\infty X \cup J_\infty C^Y$ by $\cI_X\cI_{C^Y} \subset \cI_Y$, we obtain that $(\pi_m^Y)^{-1}(T_m^o) \subset J_\infty X \setminus J_\infty W$. Of course $(\pi_m^Y)^{-1}(T_m^o) \subset S_\infty \cap (\ord_{\cJ'_Y})^{-1}(e)$, so $(\pi_m^Y)^{-1}(T_m^o) \subset T_\infty^o$, which implies that $T_\infty^o=(\pi_m^Y)^{-1}(T_m^o)$. There exists $l \ge m$ such that $\pi_m^Y(J_\infty Y)=\pi_{ml}^Y(J_lY)$ by \cite{G66}. Then $T_m^o=\pi_{ml}^Y(U_l^o)$, whence $T_m^o$ is constructible. Finally the morphism $T_{m+1}^o \to T_m^o$ is piecewise trivial with fibres $\bA^d$ by \cite[Lemma 4.1]{DL99}.
\end{proof}

We will estimate the dimension of $T_n^o$. Fix an arc $\gamma \in T_\infty^o$ such that $T_n^o$ is locally closed at $\pi_n^X(\gamma)$, and set $\gamma_m:=\pi_m^X(\gamma)$.

\begin{lemma}\label{lem:dimU}
$\dim_{\gamma_m} U_m^o \ge \dim S_m^o - c(m+1-p)$ for $m \ge n$.
\end{lemma}

\begin{proof}
Take functions $h_i$ on $A$ which define $H_i^Y$. Then a jet $\beta_m \colon \Spec k[t]/(t^{m+1}) \to A$ is contained in $J_mY$ if $\beta_m^*h_i \equiv 0$ mod $t^{m+1}$ for every $i$, but on $S_m$ we have $\beta_m^*h_i \equiv 0$ mod $t^p$ already. Hence $S_m \cap J_mY$ is obtained by cutting $S_m$ out with $c(m+1-p)$ hyperplane sections, so the lemma follows.
\end{proof}

\begin{proposition}\label{prp:dimfibre}
Let $Y$ be a l.c.i.\ scheme of dimension $d$ and $\gamma$ an arc on $Y$. Set $e:= \ord_{\cJ'_Y}(\gamma)$. Then $(\pi_{nm}^Y)^{-1}(\pi_n^Y(\gamma)) \cong \bA^{(m-n)d+e}$ for $m \ge n+e \ge 2e$.
\end{proposition}

\begin{proof}
The statement is local so we fix an embedding of $Y$ into an affine space $\bA^{d+c}=\Spec k[x_1, \ldots, x_{d+c}]$ with origin $\pi_0^Y(\gamma)$. Then $Y$ is given by $c$ functions $h_1,\ldots,h_c \in k[\underline{x}]$. By the theory of elementary divisors, there exist $P \in \GL(c,k[[t]])$ and $Q \in \GL(d+c,k[[t]])$ such that
\begin{align*}
M(\gamma)=
\begin{pmatrix}
t^{e_1} &        &         & \\
        & \ddots &         & \\
        &        & t^{e_c} &
\end{pmatrix}
\ \textrm{for} \ M=(m_{ij}):=P
\begin{pmatrix}
\frac{\partial h_1}{\partial x_1} & \cdots & \frac{\partial h_1}{\partial x_{d+c}} \\
                                  & \vdots & \\
\frac{\partial h_c}{\partial x_1} & \cdots & \frac{\partial h_c}{\partial x_{d+c}}
\end{pmatrix}
Q,
\end{align*}
where we see $\gamma$ as an element of $k[[t]]^{d+c}$ via the coordinates $x_1, \ldots, x_{d+c}$. Note that $e=\sum_i e_i$. Set ${}^t(g_1, \ldots, g_c):=P\cdot{}^t(h_1, \ldots, h_c)$, then $J_mY=\bigcap_{1 \le i \le c} (\ord_{g_i})^{-1}(\bZ_{\ge m+1})$. We also see $\gamma_n:=\pi_n^Y(\gamma)$ as an element of $(k[t]/(t^{n+1}))^{d+c} \hookrightarrow k[[t]]^{d+c}$. Consider a $(d+c)$-uple $\underline{v}={}^t(v_1, \ldots, v_{d+c})$ of polynomials in $t$ of degree less than $m-n$. Let $\underline{v}^k={}^t(v^k_1, \ldots, v^k_{d+c})$ denote the part of $\underline{v}$ of degree less than $k$ and set $\underline{v}=\underline{v}^k+t^k\underline{w}^k$. By Taylor expansion,
\begin{align*}
g_i(\gamma_n+t^{(n+1)}Q \underline{v}) =& g_i(\gamma_n+t^{(n+1)}Q \underline{v}^k) + \\
& t^{n+k+1}\sum_{1\le j \le d+c} m_{ij}(\gamma_n+t^{(n+1)}Q \underline{v}^k)w^k_j+t^{2(n+k+1)}(\cdots),
\end{align*}
whence
\begin{align*}
g_i(\gamma_n+t^{(n+1)}Q \underline{v}) \equiv g_i(\gamma_n+t^{(n+1)}Q \underline{v}^k)+t^{n+e_i+k+1}w^k_i \mod t^{n+e_i+k+2}.
\end{align*}
Therefore the condition of $\underline{v}$ to satisfy that $\gamma_n+t^{(n+1)}Q \underline{v} \in (\pi_{nm}^Y)^{-1}(\gamma_n)$ is exactly
\begin{align*}
g_i(\gamma_n+t^{(n+1)}Q \underline{v}^k)+t^{n+e_i+k+1}w^k_i \equiv 0 \mod t^{n+e_i+k+2}
\end{align*}
for $1 \le i \le c$, $0 \le k \le m-n-e_i-1$. Hence by $(\pi_{nm}^Y)^{-1}(\gamma_n) \neq \emptyset$ the proposition follows.
\end{proof}

Take $m \ge n+e$ such that $\pi_n^Y(J_\infty Y) = \pi_{nm}^Y(J_mY)$ following \cite{G66}. Then $T_n^o=\pi^Y_{nm}(U_m^o)$. Hence
\begin{align}\label{eqn:dimT}
\dim_{\gamma_n} T_n^o &= \dim_{\gamma_n} \pi^Y_{nm}(U_m^o) \\
\nonumber             &\ge \dim_{\gamma_m} U_m^o - \dim_{\gamma_m} (\pi^Y_{nm})^{-1}(\gamma_n) \\
\nonumber             &\ge \dim S_m^o - c(m+1-p) - ((m-n)d+e),
\end{align}
where the last inequality is due to Lemma \ref{lem:dimU} and Proposition \ref{prp:dimfibre}. By the inequalities (\ref{eqn:dimS}) and (\ref{eqn:dimT}) we obtain that
\begin{align*}
\dim T_n^o-(n+1)d+e > -a.
\end{align*}
By (\ref{eqn:descriptJ}) $\ord_{\cJ'_Y}=\ord_{\cJ_X}+\ord_{D^Y}$ on $J_\infty X$. Therefore Theorem \ref{thm:mld}, Lemma \ref{lem:stable} and the above inequality imply that $\mld_Z(X,D^Y)<a$, which completes the proof of Theorem \ref{thm:main}.

\section*{Acknowledgements}
I was motivated to write this paper by a question raised during the workshop at American Institute of Mathematics. I thank Doctors Shunsuke Takagi and Karl Schwede for the discussions at this workshop. American Institute of Mathematics provided financial support for my participation. This research was supported by Grant-in-Aid for Young Scientists (B) 17740015.

\bibliographystyle{amsplain}

\end{document}